\documentclass[preprint,12pt]{elsarticle}
\journal{Mathematical Biosciences}
\usepackage{amsmath,graphicx, amssymb,color}
\newtheorem{theorem}{Theorem}

\newtheorem{definition}[theorem]{Definition}

\newcommand{\R}{\mathbb{R}}
\newcommand{\M}{\mathcal{M}}
\newcommand{\LL}{\mathcal{L}}
\newcommand{\T}{\mathbb{T}}

\begin{document}
\begin{frontmatter}

\title{Pathogen evolution in
switching environments: a hybrid dynamical system approach}

\author{J\'{o}zsef Z. Farkas}
\address{Division of Computing Science and Mathematics,
University of Stirling, Stirling, FK9 4LA, United Kingdom; 
\texttt{jozsef.farkas@stir.ac.uk} }

\author{Peter Hinow}
\address{Department of
Mathematical Sciences, University of Wisconsin -- Milwaukee, P.O.~Box 413,
Milwaukee, WI 53201, USA; 
\texttt{hinow@uwm.edu} }

\author{Jan Engelst\"{a}dter}
\address{Institute of Integrative Biology,
Swiss Federal Institute of Technology, Universit\"{a}tsstrasse 16,
8092 Z\"{u}rich, Switzerland; \texttt{j.engelstaedter@uq.edu.au}\footnote{Current address: School of Biological Sciences, University of Queensland, Brisbane QLD 4072, Australia}}

%\subjclass{60K37,  92D15}
\begin{keyword}
Hybrid switching system;  pathogen evolution; stability in probability.
\end{keyword}
\date{\today}

\begin{abstract}
\begin{sloppypar}
We propose a hybrid dynamical system approach to model the evolution of a pathogen that
experiences different selective pressures according to a stochastic process. In
every environment, the evolution of the pathogen is described by a version of the  
Fisher-Haldane-Wright equation while the switching between environments follows
a Markov jump process. We investigate how the qualitative behavior of a simple single-host deterministic system changes when the stochastic switching process is added. In particular, we study the stability in probability of monomorphic equilibria.  We prove that in a ``constantly'' fluctuating environment, the genotype with the highest mean fitness is asymptotically stable in probability while all others are unstable in probability. However, if the
probability of host switching depends on the genotype composition of the
population, polymorphism can be stably maintained. 

{\bf Remark.} This is a corrected version of the paper that appeared in {\it Mathematical Biosciences} {\bf 240} (2012), p.~70-75. A corrigendum has appeared in the same journal.
\end{sloppypar}
\end{abstract}
\end{frontmatter}

\section{Introduction}
Living organisms face  changing environmental conditions.
Parasites are a case in point: after each transmission event they find
themselves in a new host that may be quite different from the previous one. For
example, the immune system of the new host may respond differently to the
parasite, and the new host may have a different genotype or even belong to a
different species. Thus, a parasite genotype that is characterized by a high
fitness in one host may have a low fitness in a different host. The question
therefore arises how parasites evolve under the fluctuating selective pressures
imposed on them through transmission events to different hosts.

Most of
the studies so far have focused on models for host-pathogen interactions in a
deterministic context \cite{RS,TB}. In some applications however it is natural
to assume that environment (and hence fitness landscape) switching is not
deterministic. For example, a pathogen could switch to a different host.
Evolution of the pathogen then takes place in the new host (or environment),
where the pathogen genotypes face different selective pressures, hence the
dynamics of the pathogen genotypes are different. We remark that the
evolving organism need not be a pathogen, nor is the environment necessarily a
``host''.

Evolution of organisms in deterministically and randomly varying
environments has been studied by many authors, see \cite{F76}
for an early review. Karlin and collaborators
\cite{KL1,KL2} introduced both deterministic and stochastic models for the
evolution of haploid and diploid organisms under changing selection intensities
for fixed and varying population sizes. In case of a deterministic two-allele
model they showed that the genotype with higher selection intensity goes to
fixation and the time to fixation varies according to the selection intensities.
Furthermore, they investigated a stochastic model where generational selection
intensities are identically distributed independent random variables. They
focused on the question how the probabilities of fixation and the times to
fixation change in the stochastic model. In \cite{Kirzhner}, Kirzhner \textit{et
al.}~considered a 4-dimensional system of difference equations for the haplotype
frequencies of a two-locus model. Typical two-locus models show either fixation
in one or both loci or stable polymorphic cycles, with period equaling the
period of the environmental changes, i.e.~the periodic fitness values. They
however showed the existence of so called supercycles that have 1100 times the
period of the periodically changing environment. The questions of structural
stability, i.e.~sensitivity in terms of the fitness parameters and the size of
the basin of attraction of these cycles were investigated. Similarly, Nagylaki
\cite{NAGY} investigated the existence of genetic polymorphisms for two-allele
models with periodically varying fitness values. He showed that in a continuous
differential equation model genetic polymorphism will persist with periods
equaling the periods of the varying fitness values, however in a discrete model
fixation is also possible.

Hybrid switching differential equations and more generally hybrid switching
diffusions have found many applications in wireless communications, queuing
networks, ecology \cite{ZY} and financial mathematics, to name but a few; see
\cite{YinZhu} and the references therein. The word ``hybrid'' refers to the
coexistence of continuous dynamics and discrete events, see also the related
concept of piecewise deterministic processes \cite{BR}. In this paper we study
a simplified version of the continuous time Fisher-Haldane-Wright equation
(also known as \textit{standard replicator equation}, \cite{F10})
subject to fitnesses  driven by a Markov jump process. It is well known that
in the deterministic Fisher-Haldane-Wright model the pathogen genotype that has the highest fitness value will
go to fixation. The coupling of different Fisher-Haldane-Wright equations by a Markov
process however requires a new definition of the concept of ``highest fitness''.  Hence we study the possible changes to the stability behaviors
of the monomorphic equilibrium states depending on the stationary
distribution of the switching Markov process. First we establish analytical
results for the stability/instability of equilibria in the hybrid model. We show that in the case of a state-independent switching process, the monomorphic equilibrium of the genotype with highest
mean fitness is asymptotically stable in probability while the monomorphic equilibria of all other genotypes are unstable in probability. As the  stationary distribution of the
switching Markov process varies, so does the  mean fitness of each genotype.
This results in exchanges of stability without the merger of the equilibria
during the transition process. We may call this a ``stochastic transcritical
bifurcation''. Finally, we  present some numerical simulations to illustrate our result and also an example of a state-dependent switching process.

\section{The switching differential equation}
We consider a model for $m$ genotypes of a pathogen evolving in $n$
possible environments. Let $w^k_i>0$ denote the fitness value of genotype $i$ in
environment $k$. We assume for simplicity that for any fixed environment $k$,
the fitness values $w^k_i$ are all different. We write $\mathbf{w}_i$ for the vector of 
all fitness values of genotype $i$. Let $P_i$ denote the
frequency of pathogen $i$, so that the dynamics in each environment takes place
in the $(m-1)$-dimensional simplex 
\begin{equation*}\T^{m-1}=\left\{P\in\mathbb{R}_{\ge0}^m \:|\: \sum_{i=1}^m P_i
= 1\right\}.
\end{equation*}
We write $P(t)=(P_1(t),P_2(t),\dots,P_m(t))$ for the state of the system at
time $t$. The frequency dynamics of pathogen genotype $i$ in environment $k$ is
given by 
\begin{equation}\label{simplex_dynamics}
\frac{dP_i(t)}{dt}=P_i(t)
\left(w_i^k-\displaystyle\sum_{j=1}^mw^k_jP_j(t)\right)=:F_i(P(t),k).
\end{equation}
This equation is the Fisher-Haldane-Wright equation for frequencies of genotypes
of asexually proliferating organisms. The rate of growth or decay of a genotype
is determined by the difference of its fitness and the average fitness of the
population. Observe that the simplex $\T^{m-1}$ and any of its subsimplices are
invariant under the dynamics given by equation \eqref{simplex_dynamics}. It can
be shown by  straightforward computation that the average fitness in environment
$k$
\begin{equation}\label{average_fitness}
\phi^k(P) = \sum_{i=1}^mw^k_iP_i 
\end{equation}
satisfies
\begin{equation*}
\begin{aligned}
\frac{d}{dt}\phi^i(P(t))&=\sum_{i=1}^mw^k_iP_i\left(w^k_i - \sum_{j=1}^mw^k_jP_j
\right) = \sum_{i=1}^mP_i(w_i^k)^2-\left(\sum_{i=1}^mP_iw_i^k\right)^2 \\
&= \sum_{i=1}^mP_i(w_i^k)^2-2\left(\sum_{i=1}^mP_iw_i^k\right)^2 +
\sum_{i=1}^m\left(\sum_{j=1}^mP_jw_j^k\right)^2 P_i  \\
&= \sum_{i=1}^mP_i\left(w^k_i - \sum_{j=1}^mw^k_jP_j
\right)^2\ge0, 
\end{aligned}
\end{equation*}
with equality if and only if $P$ is  an equilibrium. It follows from the global
existence of solutions and LaSalle's theorem that every trajectory of
\eqref{simplex_dynamics} approaches one of the finitely many equilibria
situated at the vertices of the simplex, see \cite{LA}.

The environment switches according to a continuous time stochastic process
$\alpha(t)$ that takes values in the set $\M=\{1,2,\dots,n\}$.  The switching
process $\alpha$ is a Markov process with (possibly state-dependent) generator
matrix $Q(P)$ whose entries $q_{kl}(P)$ are defined by
\begin{equation}\label{transitions}
\mathbf{P}\{\alpha(t+\Delta t) = l\,|\, \alpha(t)=k,\, (P(s),\alpha(s)),\, s\le
t \} = q_{kl}(P(t))\Delta t + o(\Delta t).
\end{equation}
The elements $q_{kl}$ of the generator matrix $Q$ satisfy $q_{kl}\ge
0$ for all $k\neq l$ and $\displaystyle\sum_{l\in\M}q_{kl}=0$ for every
$k\in\M$ (such a matrix is said to have the \textit{$q$-property}, see
\cite{YinZhang}). The complete hybrid switching ordinary differential equation
can be cast in the form
\begin{equation}\label{switching_ode}
\begin{aligned}
\frac{dP}{dt}&= F(P(t),\alpha(t)), \\
 P(0)&= p\in \T^{m-1}, \quad \alpha(0)
= \alpha \in\M, \quad \text{a.~s.},
\end{aligned}
\end{equation}
where $\alpha(t)=k$ determines the environment $k$ at time $t$ and $F=(F_1,\dots,F_m)$ are defined by \eqref{simplex_dynamics}. The
right hand side of the differential equation in \eqref{switching_ode} is globally Lipschitz continuous on the
compact set $\T^{m-1}\times\M$. This implies global existence and uniqueness
of solutions in the sense of stochastic processes, see \cite[Theorem
2.1]{YinZhu}.

For the equation \eqref{simplex_dynamics} restricted to a fixed environment $k$,
the vertices $e_i$  of the simplex $\T^{m-1}$ (i.e.~the unit vectors of
$\mathbb{R}^m$) are all the equilibrium solutions. It is easy to show that all
but one of these equilibria are unstable and that the stable equilibrium in
environment $k$ is the one for which the fitness value $w^k_i$ is the largest.
In the following section we investigate how this result generalizes to the case that
stochastic switching is introduced. For this of course, we need to first generalize the concept of stability to switching ordinary differential equations.

\section{Stability and instability in probability}
In this section we establish  results concerning the stability and instability of monomorphic steady states of the hybrid model. We first recall the following definition
%notions of stability and instability of equilibria for hybrid
%switching systems, see \cite{KK} and 
\cite[Definition 8.1]{YinZhu}. 

\begin{definition} Let $(X^{x,\alpha}(t))_{t\ge0}$ be the solution of a hybrid
switching ordinary differential equation
\begin{equation*} 
\begin{aligned}
\dot{X}(t)&=F(X(t),\alpha(t)), \quad \\
\mathbf{P}\{\alpha(t+\Delta t) &= l\,|\, \alpha(t)=k,\, (X(s),\alpha(s)),\, s\le
t \} = q_{kl}(X(t))\Delta t + o(\Delta t),\\
X(0)&=x,\quad \alpha(0)=\alpha\;\text{a.~s.},
\end{aligned}
\end{equation*}
and let (without loss of generality) $x=0$ be an equilibrium
solution, i.e.~a solution of the equation $F(0,\alpha)=0$ for every
$\alpha\in\M$. We say that $0$ is \textbf{stable in probability} if 
\begin{equation*}
\lim_{x\to 0}\mathbf{P}\left\{\sup_{t\ge0}|X^{x,\alpha}(t)|>r\right\} = 0
\end{equation*}
for every $\alpha\in\M$ and every $r>0$. We say that $0$ is
\textbf{asymptotically stable in probability} if it is stable in probability and
\begin{equation*}
\lim_{x\to 0}\mathbf{P}\left\{\lim_{t\to\infty}X^{x,\alpha}(t)=0 \right\}=1
\end{equation*}
for every $\alpha\in\M$. Finally, $0$ is \textbf{unstable in probability} if
it is not stable in probability.
\end{definition}

For $n$-tuples of functions $g(\,\cdot\,,k)\in C^1(\R^m)$ one defines a linear operator $\LL$, 
the {\it stochastic Lie derivative}  (see \cite[Equation (8.3), p.~219]{YinZhu})
\begin{equation*}
\LL g(x,k) = F(x,k)\cdot\nabla g(x,k)+\sum_{l=1}^nq_{kl}(x)g(x,k),
\end{equation*}
where $\nabla$ denotes the gradient with respect to the $x$-variable for fixed $k\in\M$. This is a natural generalization of the derivative of a scalar function along a vector field 
well known in the theory of ordinary differential equations. The following is Proposition 8.6, \cite[p.~223]{YinZhu}.

\begin{theorem}\label{Prop86} Let $D\subset\R^m$ be a neighborhood of $0$ and assume that 
there exists a function  $V:D\times\M\to[0,\infty)$
with the following properties
\begin{itemize}
\item $V(\,\cdot\,,k)$ is continuous and vanishes only at $0$,
\item $V(\,\cdot\,,k)$ is continuously differentiable in $D\setminus\{0\}$, and
\item there exists a function $\kappa:(0,r)\to(0,\infty)$ such that for all $k\in\M$ and $|x|>\varrho$,
\begin{equation*}
\LL V(x,k) \le -\kappa(\varrho)<0.
\end{equation*}
\end{itemize}
Then the equilibrium  $x=0$ is asymptotically stable in probability.
\end{theorem}
A function that satisfies the conditions of the theorem is called a {\it Lyapunov function} (for asymptotic stability).

Throughout the remainder of  this section we consider the case of a state-independent generator matrix $Q$ with a universal stationary distribution  $\pi=(q_1,\dots,q_n)$. This is the solution of the equations
\begin{equation*}
\pi\cdot\mathbf{1} = 1, \quad \textrm{and} \quad \pi Q = \mathbf{0}.
\end{equation*}
If $q_{kl}> 0$ for all $k\neq l$ then the matrix $Q$ is {\it irreducible} and $\pi>0$ is
unique \cite[p.~21]{YinZhang}. 
\begin{theorem}\label{stab_thm} Let $P_1$ be the genotype with the highest mean fitness, that is
\begin{equation}\label{1leads}
\pi\cdot \mathbf{w}_1 >\pi \cdot\mathbf{w}_i \quad \text{for all }\,i=2,\dots,m
\end{equation}
Then the equilibrium  $e_1$ is asymptotically stable in probability.
\end{theorem}
\begin{remark}\label{partition}  For almost every stationary distribution $\pi\in\T^{n-1}$, exactly one genotype satisfies a condition similar to \eqref{1leads}.
\end{remark}
\medskip

\noindent {\bf Proof.} For $i=2,\dots,m$  we set $a^k_{i,1}=w^k_i-w^k_1$ for the difference of fitness values with respect to genotype 1 and  $\mathbf{a}_{i,1} = (a_{i,1}^1,\dots,a_{i,1}^n)$. Using the constraint $\displaystyle \sum_{j=1}^m P_j(t) =1$, we eliminate $P_1$ and obtain the {\it reduced systems}
\begin{equation*}
\frac{dP_i(t)}{dt}= a_{i,1}^kP_i(1-P_i)-P_i\sum_{j=2, \,j\neq i}^ma_{i,1}^kP_j,
\end{equation*}
for $i=2,\dots,m$ and $k=1,\dots,n$. Notice that for fixed environment $k$ the linear part of this system has a diagonal structure. We define 
\begin{equation*}
\beta_i := -\pi\cdot\mathbf{a}_{i,1}>0,
\end{equation*}
with the last inequality holding true since genotype 1 has the higher mean fitness compared to every other genotype.
For $i=2,\dots,m$ we solve the systems of equations
\begin{equation*}
Q\mathbf{c}_i=\mathbf{a}_{i,1}+\beta_i\mathbf{1}
\end{equation*}
for the vector $\mathbf{c}_i=(c_i^1,\dots,c_i^n)$ where $\mathbf{1}$ is the column vector with $n$ entries 1.  The right hand sides of these equation are orthogonal to the kernel of  $Q$ which is  spanned by $\mathbf{1}$, hence there
exist solutions.  For $i=2,\dots, m$ and  $k=1,\dots,n$, we define  
\begin{equation*}
V_i(P_i,k)=(1-\gamma c_i^k) P_i^\gamma, \quad P_i>0,
\end{equation*}
with $0<\gamma <1$ yet to be selected, in such a way that all coefficients are positive. We have
\begin{equation}\label{same_calc}
\begin{aligned}
\LL V_i(P_i,k)&=\gamma(1-\gamma c_i^k)P_i^{\gamma-1}(a_{i,1}^kP_i+o(1)) +\sum_{j=1}^nq_{kj}(1-\gamma c_i^j )P_i^\gamma  \\
&=\gamma P_i^\gamma\left((1-\gamma c_i^k)a_{i,1}^k -\sum_{j=1}^nq_{kj}c_i^j +o(1)\right)\\
&=\gamma P_i^\gamma\left((1-\gamma c_i^k)a_{i,1}^k -(a_{i,1}^k+\beta_i) +o(1)\right)\\
&=\gamma P_i^\gamma\left( -\gamma c_i^k a_{i,1}^k+\pi\cdot\mathbf{a}_{i,1} + o(1)\right),
\end{aligned}
\end{equation}
where we have made use of the fact that the row sums of $Q$ are zero. 
In order to make all the factors in parentheses negative, we have to choose \mbox{$0<\gamma<1$} such that the inequality
\begin{equation}\label{ineq_need}
\pi\cdot\mathbf{a}_{i,1}<\gamma c_i^ka_i^k
\end{equation}
holds. By assumption \eqref{1leads}, the left hand side of inequality \eqref{ineq_need}
is negative. Therefore, for those indices $i$ and $k$ for which $c_i^k a_{i,1}^k\ge0$, no condition arises for $\gamma$.
If on the other hand  $c_i^ka_{i,1}^k<0$, then we can select 
\begin{equation*}
0<\gamma < \min_{i=2,\dots m \atop k=1,\dots,n } \left\{\frac{\pi\cdot\mathbf{a}_{i,1}}{c_i^k a_{i,1}^k}\::\:c_i^k a_{i,1}^k<0\right\}.
\end{equation*}
Although the $c_i^k$ are not explicitly known, this is a minimum of finitely many positive numbers. The Lyapunov function is the sum of functions of a single variable
\begin{equation*}
V(P_2,\dots,P_m,k)=\sum_{i=2}^m V_i(P_i,k)
\end{equation*}
and the condition of Proposition 8.6 in \cite{YinZhu}  follows from the linearity of the operator $\LL$ and the choice of $\gamma$. 
 \hfill $\Box$

Instability in probability of an equilibrium can be proved  similarly. The following is Proposition 8.7, \cite[p.~223]{YinZhu}. Notice however that the Lyapunov function does not vanish but has a pole at the unstable equilibrium.

\begin{theorem}\label{Prop87} Let $D\subset\R^m$ be a neighborhood of $0$ and assume that 
there exists a function  $V:D\times\M\to[0,\infty)$
with the following properties
\begin{itemize}
\item $V(\,\cdot\,,k)$ is continuously differentiable in $D\setminus\{0\}$, and
\item there exists a function $\kappa:(0,r)\to(0,\infty)$ such that for all $k\in\M$ and $|x|>\varrho$,
\begin{equation*}
\LL V(x,k) \le -\kappa(\varrho)<0,
\end{equation*}
\item for all $k\in\M$,
\begin{equation*}
\lim_{|x|\to 0} V(x,k) =\infty.
\end{equation*}
\end{itemize}
Then the equilibrium  $x=0$ is unstable in probability.
\end{theorem}

\begin{theorem}\label{instab_thm} Under the assumption \eqref{1leads} of Theorem \ref{stab_thm}, the equilibrium  $e_i,\,i=2,\dots,m$ is unstable in probability.
\end{theorem}

\noindent {\bf Proof.} The proof is very similar to that  of Theorem \ref{stab_thm}, so we only give a sketch here. This time it is $P_i$ that is being eliminated from the system containing $P_1$ and $P_i$. This results in the reduced systems
\begin{equation*}
\frac{dP_l(t)}{dt}= a_{l,i}^kP_l(1-P_l)-P_l\sum_{j\neq i,l}^ma_{j,i}^kP_j,
\end{equation*}
for $l\neq i$ and  $a^k_{l,i}=w^k_l-w^k_i$.
For $i=2,\dots,m$ let $\mathbf{c}_i=(c_i^1,\dots,c_i^n)$ be the solution  of 
\begin{equation*}
Q\mathbf{c}_i=\mathbf{a}_{1,i}-\beta_i\mathbf{1}.
\end{equation*}
 We set
\begin{equation*}
V(P_1,\dots,P_{i-1},P_{i+1},\dots,P_m,k)= V(P_1,k)= (1-\gamma c_i^k) P_1^\gamma, \quad P_1>0,
\end{equation*}
where $0>\gamma>-1$ has yet to be selected, small enough that all coefficients are positive. With a calculation similar to \eqref{same_calc} we obtain 
\begin{equation*}
\begin{aligned}
\LL V(P_1,k) 
&=\gamma(1-\gamma c_i^k)P_1^{\gamma-1}(a_{1,i}^kP_1+o(1)) +\sum_{j=1}^nq_{kj}(1-\gamma c_i^j )P_1^\gamma  \\
&=\gamma P_1^\gamma\left((1-\gamma c_i^k)a_{1,i}^k -\sum_{j=1}^nq_{kj}c_i^j +o(1)\right)\\
&=\gamma P_1^\gamma\left((1-\gamma c_i^k)a_{1,i}^k -(a_{1,i}^k-\beta_i) +o(1)\right)\\
&=\gamma P_1^\gamma\left( -\gamma c_i^k a_{1,i}^k+\pi\cdot\mathbf{a}_{1,i} + o(1)\right).
\end{aligned}
\end{equation*}
In order to make all the factors in parentheses positive (so that the entire expression becomes negative), we need to have
\begin{equation*}
0>\gamma > \max_{i=2,\dots m \atop k=1,\dots,n } \left\{\frac{\pi\cdot\mathbf{a}_{1,i}}{c_i^k a_{1,i}^k}\::\:c_i^k a_{1,i}^k<0\right\}.
\end{equation*}
The expressions whose maximum is taken are all negative since $\pi\cdot\mathbf{a}_{1,i}>0$ by assumption \eqref{1leads}.  The condition of Proposition 8.7 in \cite{YinZhu} is thereby verified. 
\hfill $\Box$

\begin{remark}\label{local} The notion of ``highest mean fitness'' requires that the generator matrix $Q$ is independent of the state and so has a universal stationary distribution $\pi$. If $Q$ depends continuously on $P$, it is still possible to formulate the corresponding ``local'' stability results for the equilibria $e_i$ by taking $\pi$ to be a stationary distribution of $Q(e_i)$.
\end{remark}

\section{Numerical simulations and examples}
\begin{figure}[ht]
\begin{center}
\includegraphics[width=60mm]{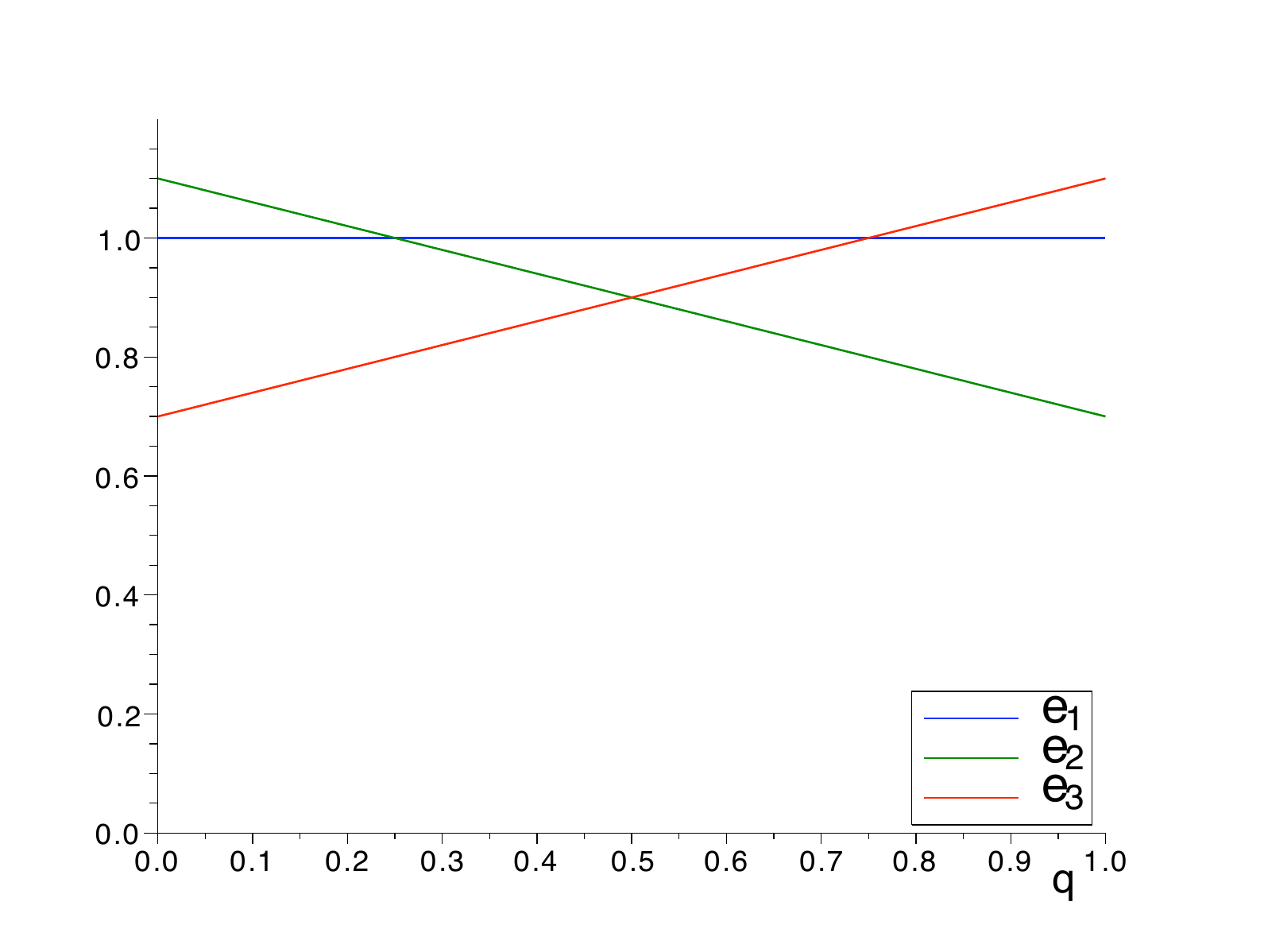}\hspace{5mm}
\includegraphics[width=60mm,height=45mm]{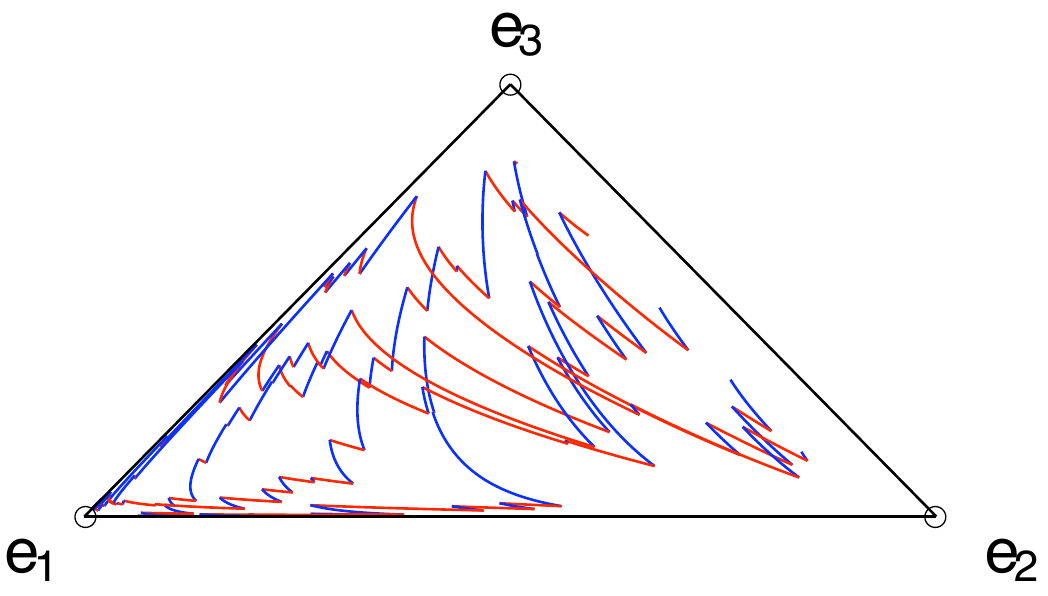}
\caption{(Left panel) The mean fitness of the three equilibria as a function of the parameter $q$ of the stationary distribution. (Right panel) Trajectories converging to $e_1$ when $q=\frac{1}{2}$. Parts of the trajectory in red indicate that environment 1 is active (when $e_3$ is attracting) while parts of the trajectory in blue indicate that environment 2 is active (when $e_2$ is attracting).}\label{saddle_stab}
\end{center}
\end{figure}

The following is an interesting example of how  stability can arise through stochastic coupling. Let the fitness values of three genotypes in two environments be given by
\begin{equation*}
\begin{aligned}
w^1_1=1,\,w^1_2=\frac{7}{10},\,w^1_3=\frac{11}{10}, \\
w^2_1=1,\,w^2_2=\frac{11}{10},\,w^2_3=\frac{7}{10}.
\end{aligned}
\end{equation*}
Although genotype 1 does not have the highest fitness in any environment, it has the highest mean fitness for stationary distributions $(q,1-q)$ with $\frac{1}{4}<q<\frac{3}{4}$, see Figure \ref{saddle_stab}. If the generator matrix $Q$ of the Markov process does not depend explicitly on the state $P(t)$, we can determine the switching times \textit{a priori} according to
$t_{i+1}=t_i+\tau$, where $\tau$ is an exponentially distributed random
variable with mean  $1$ (for example).  

\begin{figure}[ht]
\begin{center}
 \includegraphics[width=62mm]{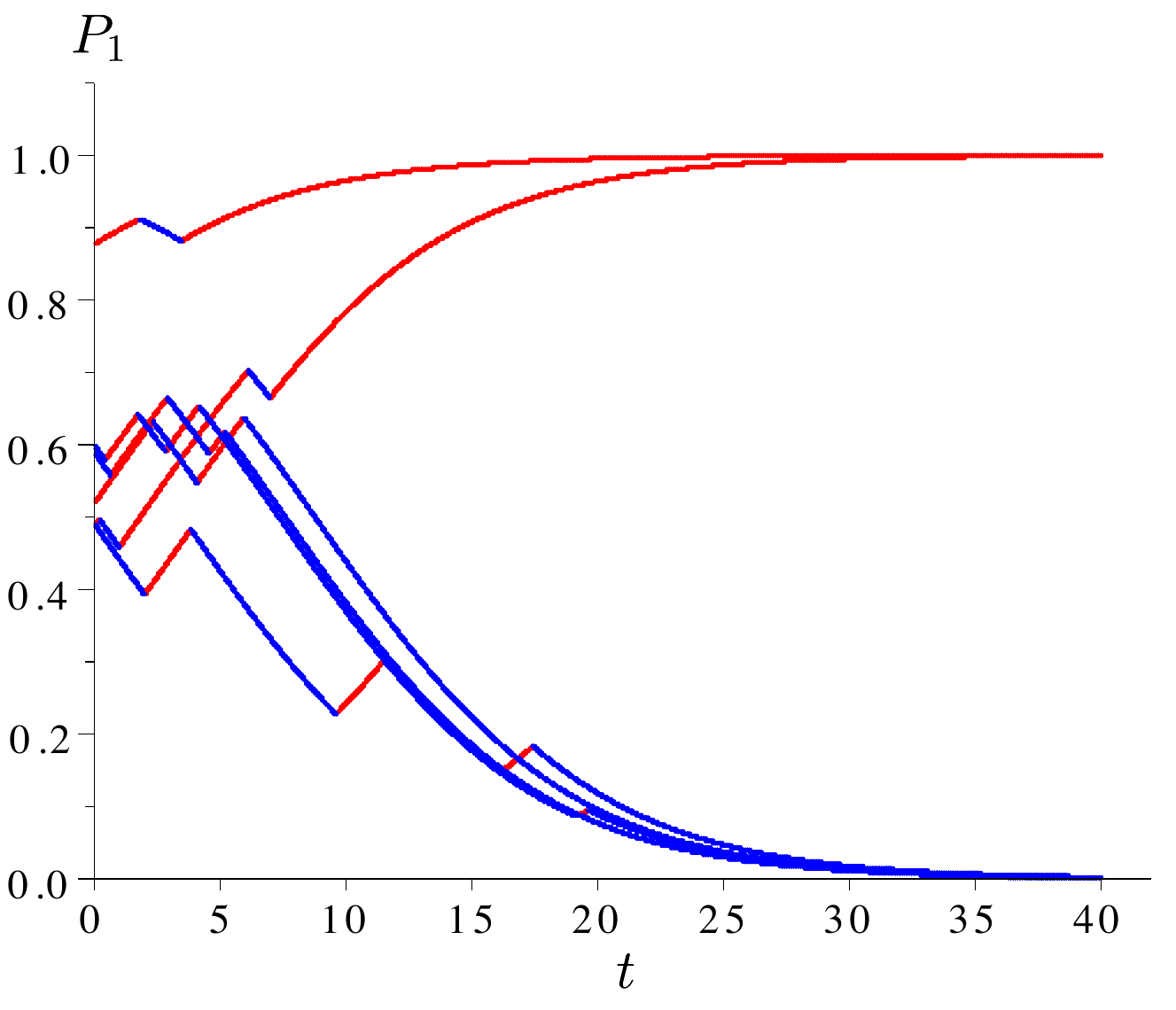}
\includegraphics[width=62mm]{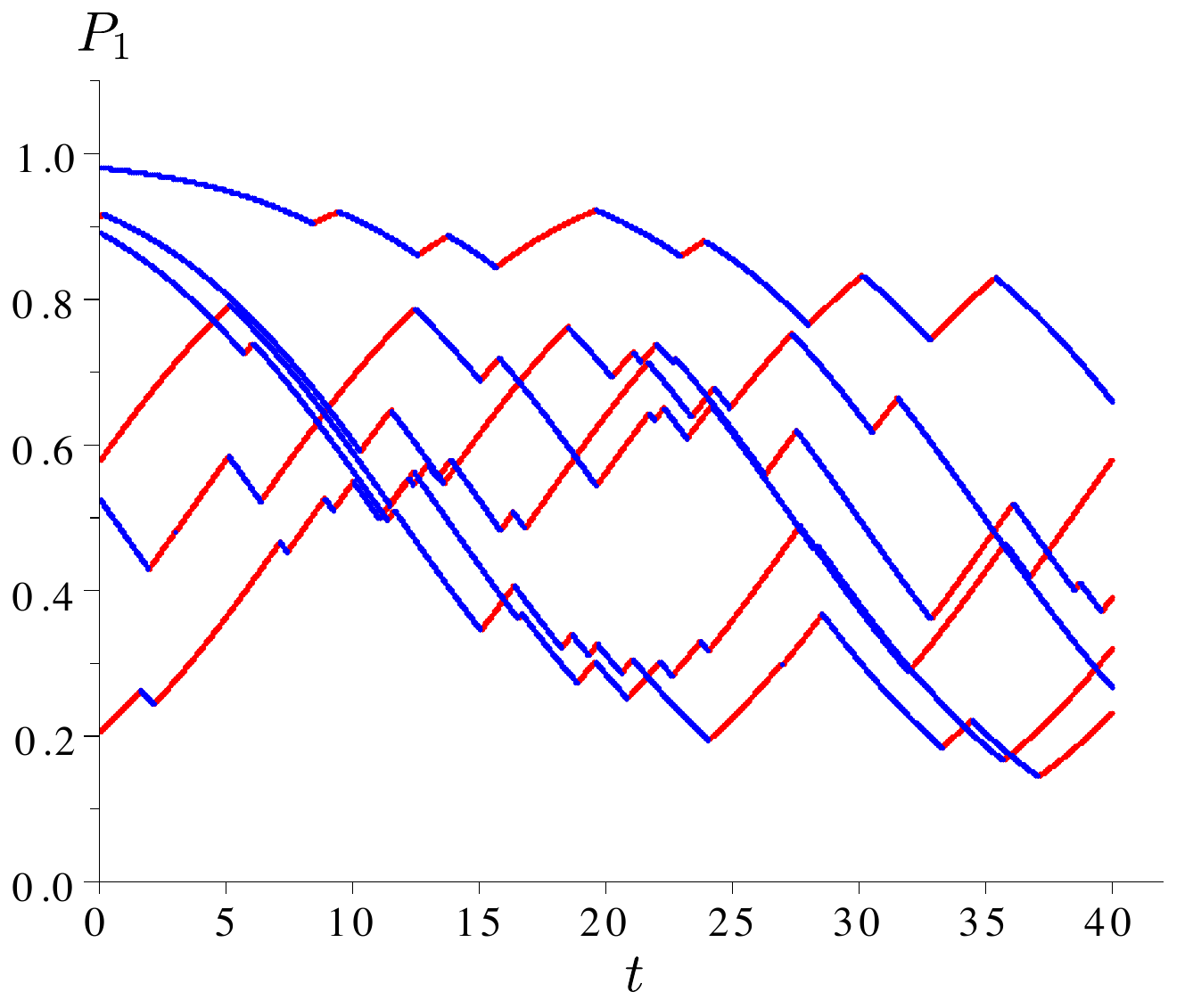}
\caption{A state-dependent generator matrix $Q_1(P)$ leads to bistability (left
panel) whereas matrix $Q_2(P)$ results in failure to converge to an
equilibrium (right panel). Parts of the trajectory in red indicate that environment 1 is active (when $e_1$ is attracting) while parts of the trajectory in blue indicate that environment 2 is active (when $e_2$ is attracting).
}\label{bistab_fig}
\end{center}
\end{figure}

To finish this section, we present two example with a state-dependent
generator matrix $Q(P)$. Let $n=m=2$,
\begin{equation*}
w^1_1 =w^2_2 =1,\:w^1_2=w^2_1=\frac{8}{10},\\
\end{equation*}
and define two switching matrix functions
\begin{equation*}
\begin{aligned}
Q_1(P_1,P_2)&=\begin{pmatrix} -P_2 & P_2 \\P_1 & -P_1 \end{pmatrix}, \\ 
Q_2(P_1,P_2)&=\begin{pmatrix} -P_1 & P_1 \\P_2 & -P_2 \end{pmatrix}.
\end{aligned}
\end{equation*}
This choice of the generator matrix means that the jump process favors jumps
into the environment that is beneficial (in case $Q_1$), respectively disadvantageous 
(in case $Q_2$) for the genotype that currently dominates.  In contrast to the previous
simulations with state-independent generator matrix, it is now necessary to
update the transition matrix of the Markov chain, namely $\exp(Q(P)\Delta t)$ during
each time step of length $\Delta t$. Following \cite[Chapter 5.3]{YinZhu}, we
use the approximation $I+Q(P)\Delta t$. 
The stationary
distribution of $Q_1(e_i)$  is, incidentally, $e_i$ for $i = 1, 2$. It follows
from Theorem \ref{stab_thm} and Remark \ref{local} that both equilibria are locally asymptotically stable in probability. Conversely, the stationary distribution of $Q_2(e_1)$ is
$e_2$ and vice versa. Under this regime,  both equilibria are locally  unstable in probability.
The results in Figure \ref{bistab_fig}
show that stochastic bistability may arise (for the choice $Q_1(P)$, left panel) or that
solutions do not converge to a monomorphic steady state (for the choice
$Q_2(P)$, right panel). 

In terms of biological interpretation, one can conceive of competing
pathogen genotypes  that cause different behaviors in the affected host.
For example, if the dominating pathogen genotype has only mild effects on
their host's well-being, infected individuals may retain their usual
mobility and thereby make a transition into a new environment more likely.
On the other hand, if the dominating pathogen genotype causes severe
morbidity, the host may exhibit restricted mobility so that a transition
into a new environment becomes less likely.

\section{Conclusions}
\begin{sloppypar}
In this work we consider the dynamics of a simple host-pathogen system, where
pathogen genotype frequencies evolve according to a simple deterministic model.
The selective pressures switch according to a Markov process. We use the
framework of switching differential equations to compare the evolution of the
pathogen in a single deterministic versus a hybrid system. In the switching
system interesting new stability patterns emerge, depending on the stationary
distribution of the underlying Markov process. We assume a fixed number of
environments (and corresponding fitnesses), in contrast to previous works. For
example, Karlin and collaborators \cite{KL1,KL2} assumed that the fitnesses
during each generation are independent identically distributed random variables.
Gillespie on the other hand in \cite{G72} proposed a stochastic differential
equation where the fitness is a process with continuous sample paths. 
\end{sloppypar}

In the case of a state-independent generator matrix of the Markov process $Q$, we have a partition of the simplex $\T^{n-1}$ of all possible stationary distributions $\pi$ into regions where one genotype has a greater mean fitness than all others, except for a set of measure zero where two genotypes have equal mean fitness (where bifurcations occur). This complete classification relies on the diagonal structure of the Jacobian of  the reduced system \eqref{reduced_systems}  at the equilibrium $0$. Due to this decoupling, it is possible to use a sum of Lyapunov functions that all depend on one variable only. In this way, we obtain a condition for asymptotic stability using convex combinations of corresponding elements of the spectra of the Jacobians in the different environments. We expect such a result to hold in the greater context of switching ordinary differential equations and diffusion processes 
with regime switching.

Our work can be refined and extended in various ways. Firstly, we use a very
simple deterministic competition model \eqref{simplex_dynamics}, where the
pathogen genotypes are ordered according to their fitness values and the only
equilibria are the vertices of the simplex $\T^{m-1}$. A straightforward
extension would be to consider the continuous time Fisher-Haldane-Wright
equation for diploid organisms for which there exist equilibria in the interior
of $\T^{m-1}$. Other competition models may lead to deterministic bistability or
to periodic
orbits (for example the rock-paper-scissors game \cite{F10}). Secondly, 
although the host switching process is stochastic, we model within-host
evolution in a deterministic way. A more realistic approach would incorporate
random genetic drift into the model. This may be particularly important during
transmission events, which often involve population bottlenecks due to small
inoculum sizes. Finally, our model only considers a single chain of transmission
events and neglects between-host selection as well as superinfections. It may be
possible to also consider multiple (branching and coalescing) transmission
chains and thus fully couple within-host and epidemiological
dynamics \cite{MAD}.

\section*{Acknowledgments}
J\'{o}zsef Z.~Farkas was partially supported by a Royal Society of Edinburgh
Grant and a University of Stirling Research and Enterprise Support Grant.  Peter
Hinow is partially supported by NSF grant DMS-1016214 and thanks the University of Stirling for its hospitality. Part of this work was
done while Jan Engelst\"{a}dter and J\'{o}zsef Z.~Farkas visited the University
of Wisconsin - Milwaukee. Financial support from the Department of Mathematical
Sciences at the the University
of Wisconsin - Milwaukee is greatly appreciated. We thank Professor Chao Zhu (University of
Wisconsin - Milwaukee) for helpful discussions 
and two  reviewers for their comments that greatly helped to improve the paper.

\end{document}